\documentclass[11pt]{article}
\usepackage{amsmath}
\usepackage{amsfonts}
\usepackage[T1]{fontenc} 
\newtheorem{theo}{Theorem}[section]
\newtheorem{lem}[theo]{Lemma}
\newtheorem{Rems}[theo]{Remarks}

\newcommand{\eq}{{\rm eq}}
\newcommand{\on}{{\rm on}}
\newcommand{\off}{{\rm off}}

\newcommand{\CQFD}{\hfill $\square$}

\newfont{\msbm}{msbm10 scaled\magstep1}
\newfont{\msbms}{msbm7 scaled\magstep1} 

\newcommand{\bbE}{\mbox{$\mbox{\msbm E}$}}

\newcommand{\bbP}{\mbox{$\mbox{\msbm P}$}}

\newcommand{\bbR}{\mbox{$\mbox{\msbm R}$}}

\newcommand{\bbsR}{\mbox{$\mbox{\msbms R}$}}

\begin{document}

\title{The on-off network traffic model under\\ intermediate scaling}
\author{Cl\'ement Dombry\footnote{Laboratoire LMA, Universit\'e de
Poitiers, T\'el\'eport 2, BP 30179, F-86962 Futuroscope-Chasseneuil cedex,
France.  Email: clement.dombry@math.univ-poitiers.fr}\; and Ingemar
Kaj\footnote{Department of Mathematics, Uppsala University, Box 480 SE
751 06 Uppsala, Sweden. Email: ikaj@math.uu.se}} \date{}

\maketitle
\begin{quote}

\abstract{ 

\noindent
The result provided in this paper helps complete a unified picture of
the scaling behavior in heavy-tailed stochastic models for
transmission of packet traffic on high-speed communication
links. Popular models include infinite source Poisson models, models
based on aggregated renewal sequences, and models built from
aggregated on-off sources. The versions of these models with finite
variance transmission rate share the following pattern: if the sources
connect at a fast rate over time the cumulative statistical
fluctuations are fractional Brownian motion, if the connection rate is
slow the traffic fluctuations are described by a stable L\'evy process,
while the limiting fluctuations for the intermediate scaling regime
are given by fractional Poisson motion.

} \end{quote}




\section{Introduction}

It is well-known that packet traffic on high-speed links exhibit data
characteristics consistent with long-range dependence and
self-similarity. To explain the possible mechanisms behind this
behavior, various network traffic models have been developed where
these features arise as heavy-tailed phenomena; see Resnick (2007) \cite{R07}. 
A natural basis for modeling such systems, applied early on during these developments, is
the view of packet traffic composed of a large number of aggregated
streams where each source alternates between an active on-state
transmitting data and an inactive off-state. The traffic streams
generate on average a given mean-rate traffic, they have stationary
increments and they are considered statistically independent. In
particular, the transmission channel is able to accommodate peak-rate
traffic corresponding to all sources being in the on-state.  To
capture in this model the strong positive dependence manifest in
empirical trace data measurements, it is assumed that the duration of
on-periods and/or off-periods are subject to heavy-tailed probability
distributions.  It is then interesting to analyze the workload of
total traffic over time and understand the random fluctuations around
its cumulative average. Our continued interest in these questions
comes from the finding that several scaling regimes exist with
disparate asymptotic limits.

The first result of the type we have in mind is Taqqu, Willinger,
Sherman (1997) \cite{TWS97}, which introduces a double limit technique. In
this sequential scheme, if the on-off model is averaged first over the
level of aggregation and then over time the resulting limit process is
fractional Brownian motion.  As the fundamental example of a Gaussian
self-similar process with long-range dependence, this limit preserves
the inherent long-range dependence of the original workload
fluctuations.  On the other hand, averaging first over time and then
over the number of traffic sources the limit process is a stable L\'evy
process. This alternative scaling limit is again self-similar but
lacks long memory since the increments are independent. Moreover,
having infinite variance the limiting workload is itself
heavy-tailed. In Mikosh, Resnick, Rootz\'en, Stegeman (2002), \cite{MRRS},
the double limits are replaced by a single scheme where instead the
number of sources grows at a rate which is relative to time. Two limit
regimes of fast growth and slow growth are identified and two limit
results corresponding to these are established, where again fractional
Brownian motion and stable L\'evy motion appear as scaled limit
processes of the centered on-off workload. The purpose of this paper
is to show that an additional limit process, fractional Poisson
motion, arises under an intermediate scheme which can be viewed as a
balanced scaling between slow and fast growth. In this case the scale
of time grows essentially as a power function of the number of traffic
sources. As will be recalled, fractional Poisson motion does indeed
provide a bridge between fractional Brownian motion and stable L\'evy
motion.

The intermediate limit scheme discussed here is indicated in Kaj
(2002) \cite{K02}, and introduced in Gaigalas and Kaj (2003)
\cite{GK}, where limit results are given for a different but related
class of traffic models under three scaling regimes referred to as
slow, intermediate and fast connection rate. The workload process is
again the superposition of independent traffic streams with
stationary increments but now each source generates packets according
to a finite mean renewal counting process with heavy-tailed
interrenewal cycle lengths. The link to the class of on-off models is
that each pair of an on-period and a successive off-period forms a
renewal cycle and the number of such on-off cycles generate a
heavy-tailed renewal counting process. Moreover, if we associate with
each renewal cycle a reward given by the length of its on-period and
apply a suitable interpretation of partial rewards, then the
corresponding renewal-reward process coincides with the on-off
workload process.

To explain briefly the limit result in \cite{GK} under intermediate
connection rate, let $(N^i(t))_i$ be i.i.d. copies of a stationary
renewal counting process associated with a sequence of inter-renewal
times of finite mean $\mu$ and a regularly varying tail function $\bar
F(t)\sim L(t)t^{-\gamma}$, characterized by an index $\gamma$,
$1<\gamma<2$, and a slowly varying function $L$.  Let $m\to\infty$ and
$a\to\infty$ in such a way that $mL(a)/a^{\gamma-1}\to \mu
c^{\gamma-1}$ for some constant $c>0$. Then the weak convergence
holds,
\[
\frac{1}{a} \sum_{i=1}^m (N^i(at)-\frac{at}{\mu})\Longrightarrow
-\frac{1}{\mu}\, c\,Y_\gamma(t/c),
\]
where $Y_\gamma(t)$ is an almost surely continuous, positively skewed,
non-Gaussian and non-stable random process, which is defined by a
particular representation of the characteristic function of its
finite-dimensional distributions.  Additional properties of the limit
process are obtained in Kaj (2005) \cite{K05} and Gaigalas (2006)
\cite{G06}, where it is shown with two different methods that
$Y_\gamma$ can be represented as a stochastic integral with respect to
a Poisson measure $N(dx,du)$ on $\bbR\times \bbR^+$ with intensity
measure $n(dx,du)= \gamma dx u^{-\gamma-1}du$. Indeed, 
\[
Y_\gamma(t)=\int_{\bbsR\times \bbsR^+} \int_0^t
            1_{[x,x+u]}(y)\,dy \,(N(dx,du)-\gamma \,dx u^{-\gamma-1}du),
\quad t\ge 0,
\]
(\cite{GK} uses $\bar F(t)\sim L(t)\gamma^{-1} t^{-\gamma}$,
consequently $n(dx,du)=dx u^{-\gamma-1}du$).  We call this process
fractional Poisson motion with Hurst index $H=(3-\gamma)/2\in
(1/2,1)$.  With 
\begin{equation}\label{eq:variancecoeff}
\sigma^2_\gamma=\frac{2}{(\gamma-1)(2-\gamma)(3-\gamma)}
=\frac{1}{2H(1-H)(2H-1)}=\sigma^2_H,
\end{equation}
we may put $Y_\gamma(t)=\sigma_\gamma P_H(t)$ and obtain 
the standard fractional Poisson motion $P_H$. A calculation reveals
\[
{\rm Cov}(P_H(s),P_H(t))=\frac{1}{2}(|s|^{2H}+|t|^{2H}-|t-s|^{2H}).
\]
For comparison, fractional Brownian motion of index $H$ has the
representation
\[
B_H(t)=\frac{1}{\sigma_H}\int_{\bbsR\times \bbsR^+} \int_0^t 
1_{[x,x+u]}(y)\,dy  \,M(dx,du),
\]
where $M(dx,du)$ is a Gaussian random measure on $\bbR\times\bbR^+$
which is characterized by the control measure $(3-2H)\,dx
u^{-2(2-H)}du$. The covariance functions of $B_H$ and $P_H$ coincide.
The fast connection rate limit for the model of aggregated renewal
processes applies if $mL(a)/a^{\gamma-1}\to \infty$ and the slow
connection rate limit if $mL(a)/a^{\gamma-1}\to 0$. For suitable
normalizing sequences, the limit processes under these assumptions are
fractional Brownian motion with Hurst index $H=(3-\gamma)/2$ in the
case of fast growth and a stable L\'evy process with self-similarity
index $1/\gamma$ in the slow growth situation, see \cite{GK}.

A number of other models have been suggested for the flow of traffic
in communication networks. The superposition of independent
renewal-reward processes applies more generally to sources which
attain random transmission rates at random times, and not merely
switch between on and off. For a model where the length of a
transmission cycle as well as the transmission rate during the cycle
are allowed to be heavy-tailed, Levy and Taqqu (2000) \cite{LT00},
Pipiras and Taqqu (2000) \cite{PT00}, and Pipiras, Taqqu and Levy
(2004) \cite{PTL04}, established results for slow and fast growth
scaling analogous to those for the on-off model. In addition, they
obtained as a fast growth scaling limit a stable, self-similar process
with stationary but not independent increments, coined the telecom
process.  A further category of models for network traffic with
long-range dependence over time starts from the assumption that
long-lived traffic sessions arrive according to a Poisson process. The
sessions carry workload which is transmitted either at fixed rate, at
a random rate throughout the session, or at a randomly varying rate
over the session length.  Such models, called infinite source Poisson
models, are widely accepted as realistic workload processes for
internet traffic. Indeed, it is natural to assume that web flows on a
non-congested backbone link are initiated according to a Poisson
process while the duration of sessions and transmission rates are
highly variable. The conditions under which slow, intermediate and
fast scaling results exist and fractional Brownian motion, fractional
Poisson motion, stable L\'evy processes and telecom processes arise in
the asymptotic limits are known in great detail for variants of the
infinite source Poisson model, see Kaj and Taqqu (2008)
\cite{KT08}. In \cite{KT08}, $Y_\gamma$ is called the intermediate
telecom process.  Mikosh and Samorodnitsky (2007) \cite{MS} consider
scaling limits for a general class of input processes, which includes
as special cases the models already mentioned as well as other
cumulative cluster-type processes. It is shown that fractional
Brownian motion is a robust limit for a variety of models under fast
growth conditions, whereas the slow growth behavior is more variable
with a number of different stable processes arising in the limit.

Our current result completes the picture for the intermediate scaling
regime, where neither of the mechanisms of fast or slow growth are
predominant. In this case, where the system workload is under
simultaneous influence of Gaussian and stable domains of attraction,
we show that the fluctuations which build up in the on-off model are
robust and again described by the fractional Poisson motion, parallel
to what is known to be valid for infinite source Poisson and
renewal-based traffic models.  In the next section 2 we introduce
properly both the on-off model and the renewal-based model to be used
as an approximation and we state the relevant background results for
these models. In section 3 we state the main result and give the
structure of the proof. Section 4 is devoted to remaining and
technical aspects of the proof.
 
\section{The on-off model and background results}

We begin by introducing the on-off model using similar notations as in
\cite{MRRS}. Let $X_{\rm on}, X_1,X_2,\dots$ be
i.i.d. non-negative random variables with distribution $F_{\on}$
representing the lengths of on-periods.  Similarly let
$Y_{\off},Y_1,Y_2,\dots$ be i.i.d. non-negative random variables with
distribution $F_{\off}$ representing the lengths of off-periods. The
$X$- and $Y$-sequences are supposed to be independent.
For any distribution function $F$ we write $\bar F=1-F$ for the right
tail. We fix two parameters,  $\alpha_\on$ and $\alpha_\off$, such that 
\begin{equation}\label{eq:indexassumption}
1<\alpha_{\on}<\alpha_{\off}<2,
\end{equation}
and assume that
\begin{equation}\label{eq:tailassumption}
\bar F_{\on}(x)=x^{-\alpha_{\on}}L_{\on}(x)\quad\mbox{and}\quad \bar
F_{\off}(x)=x^{-\alpha_{\off}}L_{\off}(x),\quad x\to\infty, 
\end{equation}
with $L_{\on}, L_{\off}$ arbitrary functions slowly varying at
infinity. Hence both distributions $F_{\on}$ and $F_{\off}$ have
finite mean values  $\mu_{\on}$ and $\mu_{\off}$ but their variances are
infinite. Assumption (\ref{eq:indexassumption}) agrees with that of
\cite{MRRS}. However, thanks to a simple symmetry
argument, we can also cover the case $\alpha_{\on}>\alpha_{\off}$.
The case $\alpha_{\on}=\alpha_{\off}$, for which the on-off process is
an alternating renewal process, falls outside of the class of
processes we are able to study within the methodology developed here.

We consider the renewal sequence generated by alternating on- and
off-periods. For the purpose of stationarity we introduce random
variables $(X_0,Y_0)$ representing the initial on- and off-periods as
follows: let $B$, $X_{\on}^{\eq}$, $Y_{\off}^{\eq}$ be independent
random variables, independent of $\{X_{\on},(X_n),Y_{\off},(Y_n)\}$,
and such that $B$ is Bernoulli with
\[
\bbP(B=1)=1-\bbP(B=0)=\mu_{\on}/\mu,
\]
and $X_{\on}^{\eq}$ and $Y_{\off}^{\eq}$ have distribution functions
\[
F_{\on}^{\eq}(x)=\frac{1}{\mu_{\on}}\int_0^x \bar F_{\on}(s)\, ds
\quad \mbox{and} \quad 
F_{\off}^{\eq}(x)=\frac{1}{\mu_{\off}}\int_0^x \bar F_{\off}(s)\,ds,
\]
respectively. Now, let
\[
X_0=B X_{\on}^{\eq}\quad \mbox{and}\quad
Y_0=B Y_{\off}+(1-B)Y_{\off}^{\eq}.
\]
Note that $X_0$ and $Y_0$ are conditionally independent given $B$ but
not independent.  At time $t=0$ the system starts in the on-state
if $B=1$ and in the off-state if $B=0$. With this initial
distribution, the alternating renewal sequence is stationary and the
probability that the system is in the on-state at any time $t$ is
$\mu_{\on}/\mu$. Renewal events occur at the start of each on-period.
Inter-renewal times are given by the independent sequence
$Z_i=X_i+Y_i, i\geq 0$, where $Z_i$ has distribution $F=F_{\on}\ast
F_{\off}$ and mean $\mu=\mu_{\on}+\mu_{\off}$ for $i\geq 1$, and $Z_0$
has distribution function
\[
F^{\eq}(x)=\frac{1}{\mu}\int_0^x \bar F(s)\,ds.
\]
The renewal sequence $(T_n)_{n\geq 1}$ with delay $T_0$ is defined by
\[
T_n=\sum_{i=0}^{n-1}Z_i,
\]
and we denote by $N(t)$ the associated counting process
\[
N(t)=\sum_{n\geq 0}1_{(0,t]}(T_n).
\]
Note that $N(t)$ has stationary increments and expectation
$\bbE[N(t)]=t/\mu$.  Moreover, because of \eqref{eq:indexassumption},
the tail behavior of the inter-renewal times is given by
\begin{equation}\label{eq:tailinterrenewal}
\bar F(x)\sim L_{\on}(x) x^{-\alpha_{\on}},\quad  x\to \infty, 
\end{equation}
see Asmussen \cite{A00}, Chapter IX, Corollary 1.11.  The on-off input process is
the indicator process for the on-state defined by
\[
I(t)=1_{[0,X_0)}(t)+\sum_{n\geq 0} 1_{[T_n,T_n+X_{n+1})}(t),\quad t\ge 0.
\]
The source is in the on-state if $I(t)=1$ and in the off-state if
$I(t)=0$. The input process $I(t)$ is strictly stationary with mean
\[
\bbE[I(t)]=\bbP(I(t)=1)=\mu_{\on}/\mu.
\]
The associated cumulative workload defined by
\[
W(t)=\int_0^t I(s)\,ds,\quad t\geq 0
\]
is a stationary increment process with mean $\bbE[W_t]=t\mu_{\on}/\mu$.

Let $(I^j,W^j,N^j)_{j\ge 1}$ denote i.i.d. copies of the input process
$I$, the accumulative workload process $W$, and the renewal counting
process $N$ for the stationary on-off model. For $m\ge 1$, consider a
server fed by $m$ independent on-off sources. We define the cumulative
workload of the $m$-server system as the superposition process
\[
W_m(t)=\sum_{j=1}^m W^j(t),\quad t\geq 0, \qquad m\ge 1,
\]
and the renewal-cycle counting process for $m$ aggregated traffic
sources by
\[
N_m(t)=\sum_{j=1}^m N^j(t),\quad t\ge 0,\qquad m\ge 1.
\]

In this paper, we are mainly concerned with the asymptotic properties
of the cumulative workload when the number of sources, $m$, increases
and time $t$ is rescaled by a factor $a>0$. Thus, we consider the centered
and rescaled process
\[
\frac{W_m(at)-mat\mu_{\on}/\mu}{b(a,m)}
=\frac{1}{b(a,m)}\sum_{j=1}^m \int_0^{at}
(I^j(s)-\frac{\mu_{\on}}{\mu})\,ds, \quad t\ge 0,
\]
where the renormalization $b(a,m)$ will be precised in the sequel.
The asymptotic is considered when both $m\to\infty$ and
$a\to\infty$. The relative growth of $m$ and $a$ have a major impact
on the limit. Let $a=a_m$ be the sequence governing the scaling of
time and suppose $a_m\to\infty$ as $m\to\infty$ (we will often omit
the subscript $m$). Following the notation in \cite{GK}, we consider
the following three scaling regimes:
\begin{itemize}
\item fast connection rate
\begin{equation}\label{FCR}
mL_{\on}(a)/ a^{\alpha_{\on}-1} \to\infty; \tag{FCR}
\end{equation}
\item slow connection rate
\begin{equation}\label{SCR}
m L_{\on}(a)/a^{\alpha_{\on}-1}
\to 0; \tag{SCR}
\end{equation}
\item intermediate connection rate
\begin{equation}\label{ICR}
m L_{\on}(a)/a^{\alpha_{\on}-1}
\to \mu\, c^{\alpha_{\on}-1},\quad 0<c<\infty. \tag{ICR}
\end{equation}
\end{itemize}

In \cite{MRRS}, the asymptotic behavior of the cumulative total
workload is investigated under conditions \eqref{FCR} and
\eqref{SCR}. 

\begin{theo}[\rm Mikosch {\it et al.}] 
Recall assumptions (\ref{eq:indexassumption}) and (\ref{eq:tailassumption}).

\begin{itemize}
\item Under condition \eqref{FCR} and with the normalization
$b(a,m)=(a^{3-\alpha_{\on}}L_{\on}(a)m)^{1/2}$, the
following weak convergence of processes holds in the space of
continuous functions on $\bbR^+$:
\[
\frac{W_m(at)-mat\mu_{\on}/\mu}{b(a,m)}\;
\Longrightarrow \;\sigma_{\alpha_\on}\,
\frac{\mu_\on}{\mu^{3/2}}\, B_H(t),\quad t\geq 0 
\]
where $B_H(t)$ is a standard fractional Brownian motion with index
$H=(3-\alpha_\on)/2$.

\item Under condition \eqref{SCR} and with the normalization 
\[
b(a,m)=(1/\bar F_{\on})^{\leftarrow}(am):=\inf\{x\geq 0: \bar
F_{\on}(x)\le 1/am\}, 
\]
then in the sense of convergence of finite dimensional distributions,
\[
\frac{W_m(at)-mat\mu_{\on}/\mu}{b(a,m)}\;
\stackrel{fdd}{\longrightarrow} \;
\sigma_0 \frac{\mu_{\off}}{\mu^{1+1/\alpha_{\on}}}\, 
X_{\alpha_{\on},1,1}(t),\quad t\geq 0, 
\]
where $X_{\alpha_{\on},1,1}(t)$ is a standard $\alpha_{\on}$-stable
L\'evy motion totally skewed to the right, i.e. such that
$$X_{\alpha_{\on},1,1}(1)\sim S_{\alpha_{\on}}(1,1,0),$$ and
\[
\sigma_0^{\alpha_\on}=
\frac{\Gamma(2-\alpha_\on)\cos(\pi\alpha_\on/2)}{1-\alpha_{\on}}.
\]
\end{itemize}
\end{theo}

The intermediate regime for renewal processes was investigated in
\cite{GK}. The formulation adopted here is given in \cite{K05}, and is
an immediate consequence of (\ref{eq:tailinterrenewal}).

\begin{theo}[\rm Gaigalas-Kaj]\label{theo:GK} 
Under condition \eqref{ICR} and with the normalization $b(a,m)=a$, the
following convergence of processes holds:
\[
\frac{N_m(at)-mat/\mu}{a}\;\Longrightarrow\; 
-\frac{1}{\mu}\sigma_{\alpha_\on}\, c\,P_H(t/c),
\]
where $\sigma_{\alpha_\on}$ is given in (\ref{eq:variancecoeff}) and $P_H(t)$ is
the standard fractional Poisson motion
\begin{equation}\label{eq:FPM}
P_H(t)=\frac{1}{\sigma_{\alpha_\on}} \int_{\bbsR\times \bbsR^+} \int_0^t
    1_{[x,x+u]}(y)\,dy \,(N(dx,du)-dx\, \alpha_{\rm on} u^{-\alpha_\on-1}du)
\end{equation}
with Hurst index $H=(3-\alpha_\on)/2$.

\end{theo}

\section{Intermediate limit for the on-off model}

In this section, we investigate the intermediate scaling limit for the
on-off model.  The following is our main result.
\begin{theo}\label{theo:main}
Under condition \eqref{ICR} and with the normalization $b(a,m)=a$, the
following convergence of processes holds in the space of continuous
functions on $\bbR^+$:
\[
\frac{W_m(at)-mat\mu_{\on}/\mu}{a}\;\Longrightarrow\;
\sigma_{\alpha_\on}\frac{\mu_\off}{\mu} \,c\, P_H(t/c),
\]
with $\sigma_{\alpha_\on}$ in (\ref{eq:variancecoeff}) and $P_H(t)$ 
the standard fractional Poisson motion in (\ref{eq:FPM}).
\end{theo}

\begin{Rems} 
{\rm The fractional Poisson motion is not self-similar but does have a
property of aggregate-similarity, introduced in \cite{K05}, which
allows for an interpretation of the scaling parameter $c$. Consider
for each integer $m\ge 1$ the sequence $c_m=m^{1/(\alpha_{\rm
on}-1)}$. Then
\[
c_mP_H(t/c_m)\;\stackrel{fdd}{=}\;\sum_{i=1}^m P_H^i(t),
\]
where $P_H^1,P_H^2,\dots$ are i.i.d. copies of $P_H$. Consider also
the sequence $c'_m=m^{-1/(\alpha_{\rm on}-1)}$. For any $m$,
\[
\sum_{i=1}^m c'_mP^i_H(t/c'_m)\;\stackrel{fdd}{=}\; P_H(t).
\]
Hence, by tracing the limit process in Theorem \ref{theo:main} as
$c_m\to\infty$, we recover in distribution the succession of all
aggregates $\sum_{1\le i\le m} P^i_H$, $m\ge 1$. Also, by letting
$c_m'\to 0$ we find that the limit process represents successively
smaller fractions which sum up to recover fractional Poisson motion. 

These relations explain the fact that fractional Poisson motion acts
as a bridge between the stable Levy process and fractional Brownian
motion. First, $\{c^H P_H(t/c)\}$ converges weakly to $\{B_H(t)\}$, as
$c\to\infty$. Indeed, $c_m^H P_H(t/c_m)\stackrel{fdd}{=}
\frac{1}{\sqrt{m}}\sum_{1\le i\le m} P^i_H(t)$ and the central limit Theorem yields the Gaussian limit as $m\to\infty$. The required tightness property is shown in \cite{G06}. Moreover, it is shown in \cite{G06}
that $c^{1/\alpha_{\rm on}}P_H(t/c)$ converges in distribution as
$c\to 0$ to the $\alpha_{\rm on}$-stable Levy process. To see that the
limit must be $\alpha_{\rm on}$-stable, take $d=c\cdot c_m'$ for any
$c>0$. Then
\[
c^{1/\alpha_{\rm on}}P_H(t/c)\;\stackrel{fdd}{=}\;
\frac{1}{m^{1/\alpha_{\rm on}}}\sum_{i=1}^m d^{1/\alpha_{\rm on}}
P_H^i(t/d),\quad m\ge 1,
\]
and, assuming that the rescaled process $(c^{1/\alpha_{\rm on}}P_H(t/c))_{t\geq 0} $ converge to some non-trivial limit process $L$, we must have as $c\to 0$ (and hence $d\to 0$)
\[
L(t)\;\stackrel{fdd}{=}\;
\frac{1}{m^{1/\alpha_{\rm on}}}\sum_{i=1}^m L^i(t),\quad m\ge 1.
\]
This indicates that the limit $L$ must be $\alpha_{\rm on}$-stable.}
\end{Rems}

\paragraph{\bf Heuristics of the proof of Theorem \ref{theo:main}}

To motivate that the limit process in the intermediate
connection rate limit appears naturally, we discuss a
decomposition of the centered on-off process based on its
representation as a renewal-reward model. We first note that the
single source cumulative workload has the form
\[
W(t)=X_0\wedge t +\sum_{i=1}^{N(t)}X_i -(T_{N(t)-1}+X_{N(t)}-t)_+.
\]
Similarly, focusing on off-periods rather than on-periods, we have
\[
t-W(t)=Y_0\wedge t +\sum_{i=1}^{N(t)}Y_i -(T_{N(t)}-t)\wedge Y_{N(t)}.
\]
The centered single source workload is therefore
\begin{eqnarray*}
W(t)-\frac{\mu_\on}{\mu}t
&=& -(t-W(t))+\frac{\mu_{\off}}{\mu}t\\
&=& -\mu_{\off}(N(t)-t/\mu) -\sum_{i=1}^{N(t)}(Y_i-\mu_{\off})+R(t)
\end{eqnarray*}
with
\[
R(t)=(T_{N(t)}-t)\wedge Y_{N(t)}- Y_0\wedge t.
\]
Thus, for the workload of $m$ sources,
\begin{equation}\label{eq:decomposition}
W_m(t)-\frac{\mu_\on}{\mu}mt=-\mu_\off(N_m(t)-mt/\mu)
 -\sum_{j=1}^m \sum_{i=1}^{N^j(t)}(Y^j_i-\mu_\off)+\sum_{j=1}^m R^j(t)
\end{equation}
using obvious notations. The balancing of terms under the
scaling relation \eqref{ICR}, makes it plausible that both terms
\[
\frac{1}{a}\sum_{j=1}^m
 \sum_{i=1}^{N^j(at)}(Y^j_i-\mu_\off),\quad \frac{1}{a}\sum_{j=1}^m R^j(at)
\]
vanish in the scaling limit. This suggests asymptotically,
\begin{equation}\label{eq:comparison}
\frac{W_m(at)-\mu_\on mat/\mu}{a}\sim -\,\mu_\off\frac{N_m(at)-mat/\mu}{a},
\end{equation}
and so Theorem \ref{theo:GK} would imply Theorem \ref{theo:main}. In
the next final section, we will compare rigorously the two
processes in (\ref{eq:comparison}).

\section{Proof of Theorem \ref{theo:main}}

The proof of Theorem \ref{theo:main} relies on the following three lemmas:

\begin{lem}\label{lem:remainder1}
In the scaling \eqref{ICR}, for all $t\geq 0$,
\begin{equation}\label{eq:remainder1}
\frac{1}{a}\sum_{j=1}^m  \sum_{i=1}^{N^j(at)}(Y^j_i-\mu_\off) \Longrightarrow 0.
\end{equation}
\end{lem}

\begin{lem}\label{lem:remainder2}
In the scaling \eqref{ICR}, for all $t\geq 0$,
\[
\frac{1}{a}\sum_{j=1}^m R^j(at) \Longrightarrow 0.
\]
\end{lem}

\begin{lem}\label{lem:tightness}
In the scaling \eqref{ICR}, the sequence of processes 
\[
\frac{W_m(at)-mat\mu_{\on}/\mu}{a},\quad t\geq 0, \quad m\geq 1
\]
 is tight in the space of continuous functions on $\bbR^+$.
\end{lem}

\paragraph{\bf Proof of Theorem \ref{theo:main}.} 
By Theorem \ref{theo:GK}, Lemma \ref{lem:remainder1} and Lemma
\ref{lem:remainder2}, the convergence of finite-dimensional
distributions,
\[
\frac{W_m(at)-mat\mu_{\on}/\mu}{a}\;\stackrel{fdd}{\Longrightarrow}
\;\sigma_{\alpha_\on}\frac{\mu_\off}{\mu} \,c\, P_H(t/c),
\]
is a consequence of the decomposition given in
\eqref{eq:decomposition}. By Lemma \ref{lem:tightness} the sequence is
tight in the space of continuous functions on $\bbR^+$.  Hence weak
convergence holds in the space of continuous functions and Theorem
\ref{theo:main} is proved. \hfil \CQFD

\paragraph{\bf Proof of Lemma \ref{lem:remainder1}.}
We construct an alternative representation of the random variable in
the left hand side of \eqref{eq:remainder1}. Define
\[
\widetilde N^1(at)=\inf\{k\geq 0;\ X_0^1+\sum_{i=1}^k Z^1_i\geq at\}
\]
and for $j\geq 2$
\[
\widetilde N^{j}(at)=\inf\{k\geq 0;\ X_0^{j}+\sum_{i=1}^{k} Z^1_{\widetilde
N^{j-1}(at)+i}\geq at\}.
\]
For $m\geq 1$, let $\widetilde N_m(at)=\sum_{j=1}^m \widetilde N^j(at)$.  The
random variables $\widetilde N^j(at), j\geq 1$ are i.i.d and for each
fixed $t\geq 0$
\[
\frac{1}{a}\sum_{j=1}^m
\sum_{i=1}^{N^j(at)}(Y^j_i-\mu_\off)\quad \mbox{and}\quad
\frac{1}{a}\sum_{i=1}^{\widetilde N_m(at)}(Y^1_i-\mu_\off),
\] 
have the same distribution (note that
the uni-dimensional marginal distributions are equal but not the
multidimensional distributions). This representation will enable us to
prove that under assumption \eqref{ICR}, in the space of c\'ad-l\'ag
functions on $\bbR^+$ endowed with the Skorokhod topology, we have the
convergence
\begin{equation}\label{eq:conv1}
\left(\frac{1}{a}\sum_{i=1}^{ am u }(Y^1_i-\mu_\off)\right)_{u\geq 0}
 \Longrightarrow 0.
\end{equation}
Moreover, 
\begin{equation}\label{eq:conv2}
\frac{1}{am}\tilde N_m(at)\Longrightarrow \frac{t}{\mu}.
\end{equation}
Equations \eqref{eq:conv1} and \eqref{eq:conv2} together imply  
\[
\frac{1}{a}\sum_{i=1}^{\tilde N_m(at)}(Y^1_i-\mu_\off)\Longrightarrow 0
\]
and this proves the lemma.  Thus, it remains to prove \eqref{eq:conv1} and
\eqref{eq:conv2}.

To this aim, recall that the random variables $Y^1_i, i\geq 1,$ are
i.i.d. with distribution such that the tail function $\bar F_{\off}$ is
regularly varying with index $-\alpha_{\off}$. Hence there exists a
regularly varying function $L$ such that the centered and rescaled sum
\[
\Big(\frac{1}{( am)^{1/\alpha_\off}L(am)}\sum_{i=1}^{ am u
}(Y^1_i-\mu_\off)\Big)_{u\geq 0}
\]
converges in the space of c\'ad-l\'ag functions to some
$\alpha_\off$-stable L\'evy process (see \cite{P75}, the exact form of $L$ or of the
limit process are not needed here). This implies the convergence
property \eqref{eq:conv1}, since $a>\!>  (am)^{1/\alpha_\off}L(am)$
under the scaling assumption \eqref{ICR} with $\alpha_\on<\alpha_\off$.

We now prove equation \eqref{eq:conv2}. The stationary renewal process
$N(t)$ has mean $t/\mu$ and variance given asymptotically by
\[
{\rm Var}(N(t))\sim
\sigma_{\alpha_\on}^2\frac{1}{\mu^3}t^{3-\alpha_\on}L_\on(t),\quad
t\to\infty,
\]
see \cite{GK}, Equation (30), and references therein.  Hence,
$\frac{1}{am}\widetilde N_m(at)$ has mean $t/\mu$ and variance under scaling
\eqref{ICR}, such that  
\begin{eqnarray*}
{\rm Var}\left[\frac{1}{am}\widetilde N_m(at)\right]&=&\frac{1}{a^2m}{\rm
Var}[N(at)]\\ &\sim&
\frac{a^{1-\alpha_\on}L_\on(at)}{m}\sigma_{\alpha_\on}^2
\frac{1}{\mu^3}t^{3-\alpha_{\on}}\to 0.
\end{eqnarray*} 
This shows that $\frac{1}{am}\widetilde N_m(at)$ converges in
distribution to $t/\mu$, which is \eqref{eq:conv2}. This ends the
proof of Lemma \ref{lem:remainder1} \hfill \CQFD

\paragraph{\bf Proof of Lemma \ref{lem:remainder2}.} 

Since $|R^j(t)|\leq Y_0^j+(T^j_{N^j(t)}-t)\wedge Y^j_{N^j(t)}$, it is
enough to prove that
\[
\frac{1}{a}\sum_{j=1}^m Y_0^j\Longrightarrow 0 \quad {\rm and}\quad
\frac{1}{a}\sum_{j=1}^m (T^j_{N^j_t}-t)\wedge
Y^j_{N^j(t)}\Longrightarrow 0.
\]
By stationarity, the random variables $Y^j_0$ and
$(T^j_{N^j_t}-t)\wedge Y^j_{N^j(t)}$ have the same distribution; they
represent the remaining time after $0$ and $t$, respectively, of the
first off-period after $0$, and after $t$.  Since both sums have the
same distribution, we only consider the first one. 

Using Karamata's Theorem (see \cite{BGT}), the tail function $\bar F_\off^\eq$
satisfies
\[
\bar F_\off^\eq(x)=\frac{1}{\mu_\off}\int_x^{\infty} \bar F_\off(s)ds \sim \frac{1}{\mu_\off}\frac{x^{-(\alpha_\off-1)}}{\alpha_\off-1}L_\off(x)
\]
as $x\to \infty$. This implies that the random variable
$Y_0$ has a regularly varying tail with index $-(\alpha_{\off}-1)$ and
hence belongs to the domain of attraction of an $(\alpha_{\off}-1)$-stable
distribution. Therefore there exists a slowly varying function $L$, such that
\[
\frac{1}{m^{1/(\alpha_{\off}-1)}L(m)}\sum_{j=1}^m Y_0^j
\]
converges in distribution to a stable law of index
$\alpha_{\off}-1$. Under scaling \eqref{ICR} with
$\alpha_{\on}<\alpha_{\off}$, we have
$a>\!>m^{1/(\alpha_{\off}-1)}L(m)$ and hence
\[
\frac{1}{a}\sum_{j=1}^m Y_0^j\Longrightarrow 0.
\]
\vskip -12mm\hfill\CQFD

\paragraph{\bf Proof of Lemma \ref{lem:tightness}.} 
The proof given in \cite{MRRS} for fast scaling \eqref{FCR} can be
adapted to our settings. We recall only the main lines. According to
Billingsley \cite{Bill}, Theorem 12.3, it is enough to prove that for
any $t_1$, $t_2$ with $|t_1-t_2|\leq 1$ and for some $\varepsilon>0$,
there exists a constant $C>0$ and an $a_0>0$, such that for all $a\geq a_0$
\[
\bbE\left[\frac{1}{a}\left| (W_{m}(at_2)-mat_2\mu_{\on}/\mu)-
(W_{m}(at_1)-mat_1\mu_{\on}/\mu)\right|^2\right]\leq
C|t_2-t_1|^{1+\varepsilon}.
\]
Using the definition of $W_m$, centering and stationarity, it is enough to prove
that for all $t\in [0,1]$ and $a\geq a_0$,
\begin{equation}\label{eq:tightnessbound}
\frac{m}{a^2}{\rm Var}[ W_{at}]\leq Ct^{1+\varepsilon}
\end{equation}
(the constant $C$ may change from one appearance to
another).  However, according to \cite{MRRS}, Equation (7.1),
\begin{equation}\label{eq:estimW}
{\rm Var}(W_t)\sim \sigma_{\alpha_\on}^2\frac{\mu_{\on}^2}{\mu^3} 
t^{3-\alpha_{\on}}L_{\on}(t),\quad t\to\infty.
\end{equation} 
This relation 
and the scaling (ICR) together imply,
as $a\to\infty$,
\[
\frac{a^2}{m}\sim \frac{c^{1-\alpha_{\on}}}{\mu} a^{3-\alpha_{\on}}L_{\on}(a) \sim \frac{c^{1-\alpha_{\on}}}{\sigma^2_{\alpha_\on}}\frac{\mu^2}{\mu_\on^2} \ {\rm Var}[ W_{a}],  
\]
and so there is $C>0$, such that for $a$ large enough
\[
\frac{m}{a^2}{\rm Var}[ W_{at}]\leq C\frac{ {\rm Var}[ W_{at}]}{{\rm Var}[ W_{a}]}.
\]
By \eqref{eq:estimW}, the function $a\mapsto {\rm Var}[ W_{a}]$ is
regularly varying with index $3-\alpha_\on$. Then, using Potter
bounds (see \cite{BGT}), we conclude that there exist $a_0>0$ and
$\varepsilon<1-\alpha_\on/2$, such that for all $t\in
(0,1)$ and $a\geq a_0/t$,
\[
\frac{ {\rm Var}[ W_{at}]}{{\rm Var}[ W_{a}]}\leq 
\frac1{1-\varepsilon}t^{3-\alpha_\on-\varepsilon}.
\]
(see the proof of Lemma 13 in \cite{MRRS} for details). This implies that for all $t\in (0,1)$
and all $a$ such that $at\geq a_0$,
\[
\frac{m}{a^2}{\rm Var}[ W_{at}]\leq \frac{C}{1-\varepsilon}t^{3-\alpha_\on-\varepsilon}\le C t^{1+\varepsilon}.
\]
On the other hand, if $t\leq a_0/a$, then, for $a$ large enough,  
\[
\frac{m}{a^2}{\rm Var}[ W_{at}]\le \frac{Ca^2t^2}{{\rm Var}[ W_{a}]}
\le C \frac{a^2t^2}{a^{3-\alpha_\on}L_{\on}(a)}
\le C \frac{(at)^{1+\varepsilon}a_0^{1-\varepsilon}}
{a^{3-\alpha_\on}L_{\on}(a)}
\]
and so 
\[
\frac{m}{a^2}{\rm Var}[ W_{at}] \le C
\frac{a_0^{1-\varepsilon}t^{1+\varepsilon} }
{ a^{2-\alpha_\on-\varepsilon}L_\on(a) } 
 \le C\ t^{1+\varepsilon}.
\]
In the last inequality, we use the fact that  $2-\alpha_\on-\varepsilon>0$ and so $a^{2-\alpha_\on-\varepsilon}L_\on(a)\to \infty$ as $a\to\infty$; taking $a_0$ large enough, we can suppose that for $a\geq a_0$, $a^{2-\alpha_\on-\varepsilon}L_\on(a)$ remains bounded away from zero. 

By combining the estimates for the cases $at\ge a_0$ and $at\le a_0$ we obtain 
(\ref{eq:tightnessbound}), which completes the proof. \CQFD

\end{document}